\documentclass[reqno,11pt]{amsart}

\usepackage{amsmath,latexsym,amssymb}

\def\sqr#1#2{{\vcenter{\hrule height.#2pt
        \hbox{\vrule width.#2pt height#1pt \kern#1pt
                \vrule width.#2pt}
        \hrule height.#2pt}}}
\def\square{\mathchoice\sqr64\sqr64\sqr{4}3\sqr{3}3}
\def\QED{\hfill$\square$}
\def\demo{\noindent{\bf Proof: }}
\def\manyby{\hbox to.65in{\hrulefill}}

\def\lar{\longrightarrow}
\def\rar{\rightarrow}

\def\tratto{\mbox{\rule{2.5mm}{.3mm}$\;\!$}}

\newtheorem{Theorem}{\bf Theorem}[section]

\newtheorem{Corollary}[Theorem]{\bf Corollary}
\newtheorem{Proposition}[Theorem]{\bf Proposition}

\newtheorem{Remark}[Theorem]{\bf Remark}
\newtheorem{Example}[Theorem]{\bf Example}

\newtheorem{Question}[Theorem]{\bf Question}

\begin{document}

\baselineskip=14pt

\title[Sally modules of ${\mathfrak m}$-primary
ideals in local rings]{}

\author[Alberto Corso]{}

\maketitle

\pagestyle{myheadings}

\ \vspace{0.9in}

\noindent \text{\Large\bf Sally modules of ${\mathfrak m}$-primary
ideals in local rings}\footnote{AMS 2000 {\em Mathematics Subject
Classification}. Primary 13A30, 13B21, 13D40; Secondary 13H10,
13H15. \newline The author would like to thank W. Heinzer, C.
Polini and W.V. Vasconcelos for fruitful conversations concerning
the material in this article, and the University of Notre Dame for
the pleasant hospitality. }

\vspace{.25in}

\noindent {\large\sc Alberto \ Corso} \newline Department of
Mathematics, University of Kentucky, Lexington, Kentucky 40506 ---
USA \newline {\it URL}: {\tt
http:\!/\!/www.ms.uky.edu{/}{\textasciitilde}corso}
\newline {\it E-mail}: {\tt corso@ms.uky.edu}

\vspace{.4in}

\baselineskip=8pt

\begin{quote}
{\tiny {\bf Abstract:} \ Given a local Noetherian ring $(R,
{\mathfrak m})$ of dimension $d>0$ and infinite residue field, we
study the invariants $($dimension and multiplicity$)$ of the Sally
module $S_J(I)$ of any ${\mathfrak m}$-primary ideal $I$ with
respect to a minimal reduction $J$. As a by-product we obtain an
estimate for the Hilbert coefficients of ${\mathfrak m}$ that
generalizes a bound established by J. Elias and G. Valla in a
local Cohen-Macaulay setting. We also find sharp estimates for the
multiplicity of the special fiber ring ${\mathcal F}(I)$, which
recover previous bounds established by C. Polini, W.V. Vasconcelos
and the author in the local Cohen-Macaulay case. Great attention
is also paid to Sally modules in local Buchsbaum rings.}
\end{quote}

\baselineskip=14pt

\vspace{.4in}

\section{\bf Introduction \hfill\break}

\noindent Let $(R, {\mathfrak m})$ be a local Noetherian ring of
dimension $d>0$ and with infinite residue field and let $I$ be an
${\mathfrak m}$-primary ideal. The {\it Rees algebra} ${\mathcal
R}(I)$ $($often denoted $R[It]$, where $t$ is an indeterminate
over $R)$, the {\it associated graded ring} ${\mathcal G}(I)$
$($often denoted ${\rm gr}_I(R))$, and the {\it special fiber
ring} ${\mathcal F}(I)$ of $I$
\[
{\mathcal R}(I) = \bigoplus_{n=0}^{\infty} I^nt^n, \qquad
{\mathcal G}(I) = {\mathcal R} \otimes R/I, \qquad {\mathcal F}(I)
= {\mathcal R} \otimes R/{\mathfrak m},
\]
collectively referred to as {\it blowup algebras of $I$}, play an
important role in the process of blowing up the variety ${\rm
Spec}(R)$ along the subvariety $V(I)$. In particular, their depth
properties have been under much scrutiny in the past two decades.
Moreover, these algebras are also extensively used as the means to
examine diverse properties of the ideal $I$.

A successful approach to the study of blowup algebras, initiated
by S. Goto, S. Huckaba, C. Huneke, D. Johnston, D. Katz, K.
Nishida, N.V. Trung and others, uses {\it minimal reductions} of
the ideal. This notion was first introduced and exploited by D.G.
Northcott and D. Rees fifty years ago for its effectiveness in
studying multiplicities in Noetherian local rings \cite{NR}. We
recall that, in our setting, a minimal reduction $J$ of $I$ is a
$d$-generated subideal of $I$ such that $I^{r+1}=JI^r$ for some
non-negative integer $r$. Phrased otherwise, we say that any such
$J$ is a minimal reduction of $I$ if the inclusion of Rees
algebras ${\mathcal R}(J) \hookrightarrow {\mathcal R}(I)$ is
module finite and $r$ is the bound of the degrees required to
generate ${\mathcal R}(I)$ as a module over ${\mathcal R}(J)$. The
underlying philosophy is that it is reasonable to expect to
recover some of the properties of ${\mathcal R}(I)$ from the
amenable structure of ${\mathcal R}(J)$, especially whenever $r$
is sufficiently small.

An earlier approach to the depth properties of the blowup algebras
dates back to J. Sally and involves a detailed analysis of
numerical information encoded in the Hilbert-Samuel function of
$I$, that is the function that measures the growth of the length
$\lambda(R/I^n)$ of $R/I^n$ for all $n \geq 1$. It is well known
that for $n \gg 0$ the function $\lambda(R/I^n)$ is a polynomial
in $n$ of degree $d$, say
\[
\lambda(R/I^n) = e_0(I) {{n+d-1}\choose d} - e_1(I)
{{n+d-2}\choose d-1} + \cdots + (-1)^d e_d(I),
\]
where $e_0(I), e_1(I), \ldots, e_d(I)$ are called the {\it Hilbert
coefficients} of $I$.

From a more recent vintage is a remarkable novelty which bridges
the previous two approaches. More precisely, in \cite{Vas} W.V.
Vasconcelos enlarged the list of blowup algebras by introducing
the {\it Sally module} $S_J(I)$ of $I$ with respect to a minimal
reduction $J$. This is the graded ${\mathcal R}(J)$-module defined
in terms of the short exact sequence
\[
0 \rar I {\mathcal R}(J) \lar I {\mathcal R}(I) \lar S_J(I) =
\bigoplus_{n=2}^\infty I^n/IJ^{n-1} \rar 0.
\]
When $R$ is a {\it local Cohen-Macaulay} ring, he shows that the
Sally module $S_J(I)$ is a $d$-dimensional graded ${\mathcal
R}(J)$-module with a unique associated prime, namely ${\mathfrak
m}{\mathcal R}(J)$, provided $S_J(I)$ is not the trivial module.
He also finds precise relations among the Hilbert coefficients of
$I$ and $S_J(I)$, which in turn enable him to recover the bound
$e_1(I)-e_0(I)+\lambda(R/I) \geq 0$, originally due to D.G.
Northcott \cite{No}. At the same time, he obtains the result of C.
Huneke \cite{H2} and A. Ooishi \cite{O}, which says that equality
holds in Northcott's estimate if and only if $I^2=JI$ for some
(equivalently, any) minimal reduction $J$ of $I$. In particular,
it follows that ${\mathcal G}(I)$ is Cohen-Macaulay when equality
holds. Later, the Sally module has been further studied by M. Vaz
Pinto \cite{Vaz}, H.-J. Wang \cite{W1,W2,W3}, the author, C.
Polini and M. Vaz Pinto \cite{CPVP}, L. Doering and M. Vaz Pinto
\cite{DVP}, and C. Polini \cite{P}. More recently, in a joint work
with C. Polini and W.V. Vasconcelos \cite{CPV2}, we have used the
Sally module as the means to obtain information on the
multiplicity of the special fiber ring ${\mathcal F}(I)$, on the
unmixedness $($or even Cohen-Macaulayness$)$ of ${\mathcal F}(I)$,
and, ultimately, on the reduction number $r$ of $I$.

Our investigation has been prompted by the lack of knowledge of
the properties of Sally modules in non Cohen-Macaulay settings. In
this paper we study the invariants $($dimension and
multiplicity$)$ of the Sally module $S_J({\mathfrak m})$ in an
arbitrary local Noetherian ring $R$ of dimension $d>0$ and
infinite residue field. In particular, in
Proposition~\ref{dimension} we show that the Sally module $S_J(I)$
has dimension $d$ if and only if $e_1({\mathfrak
m})-e_0({\mathfrak m})-e_1(J)+1$ is strictly positive. However,
the dimension of the Sally module $S_J({\mathfrak m})$ may assume
intermediate values, as shown in Example~\ref{ex-1}. Interestingly
enough, the situation becomes rather extremal in the case of a
local Buchsbaum ring, as we show in Proposition~\ref{dimbuch} that
in this setting the Sally module $S_J({\mathfrak m})$ has either
dimension $d$ or $0$. As a consequence of
Proposition~\ref{dimension} we prove in Theorem~\ref{E-V} that in
an arbitrary local Noetherian ring $R$ the estimate
$2e_0({\mathfrak m})-e_1({\mathfrak m})+e_1(J) \leq \mu({\mathfrak
m})-d+2$ always holds. This bound was first obtained by J. Elias
and G. Valla in the case of a local Cohen-Macaulay ring \cite{EV}.

As far as the case of an arbitrary ${\mathfrak m}$-primary ideal
is concerned, we do not find a closed formula for the multiplicity
of its Sally module with respect to a minimal reduction, but we
give an upper bound. In Proposition~\ref{sally}, we show that the
multiplicity of a $d$-dimensional Sally module $S_J(I)$ is {\it at
most} $e_1(I)-e_0(I)-e_1(J)+\lambda(R/I)$, with equality if and
only if $I$ contains the zero-th local cohomology $H_{\mathfrak
m}^0(R)$ of $R$ with support in ${\mathfrak m}$. Nevertheless,
this estimate allows us to find sharp bounds for the multiplicity
$f_0(I)$ of the special fiber ring ${\mathcal F}(I)$. To be more
specific, in Theorem~\ref{fiber} we show that $f_0(I) \leq
e_1(I)-e_0(I)-e_1(J)+\lambda(R/I)+\mu(I)-d+1$, thus generalizing a
previous result obtained jointly with C. Polini and W.V.
Vasconcelos \cite{CPV2}. If in addition $R$ is a local Buchsbaum
ring, in Theorem~\ref{fiber2} we show that $f_0(I) \leq e_1(I)+
I(R)-e_1(J)+1$, where $I(R)$ is the Buchsbaum invariant of $R$
introduced by J. St\"uckrad and W. Vogel \cite{SV}. This result
generalizes a previous estimate due to W.V. Vasconcelos
\cite{mrn}. We also show that if equality holds in the latter
estimate then the ideal $I$ has minimal multiplicity in the sense
of S. Goto (see \cite{CPV2} for a similar statement). Still in a
Buchsbaum setting, we show in Proposition~\ref{dimbuch-primary}
that the Sally module of an ideal $I$ containing $H_{\mathfrak
m}^0(R)$ has either dimension $d$ or $0$.

\bigskip

\section{ \bf The case of the maximal ideal and applications
to Hilbert coefficients \hfill\break}

\noindent We first do some calculations in order to compute the
dimension and the multiplicity of the Sally module $S_J({\mathfrak
m})$ in an arbitrary local Noetherian ring $R$ of dimension $d>0$
with infinite residue field. We show in
Proposition~\ref{dimension} that the dimension is exactly $d$ if
and only if its multiplicity has a precise value, namely
$e_1({\mathfrak m})-e_0({\mathfrak m})-e_1(J)+1$. This allows us
to obtain a general estimate relating the first two Hilbert
coefficients of ${\mathfrak m}$. However, the situation becomes
more interesting in the case of a local Buchsbaum ring, as we show
in Proposition~\ref{dimbuch} that in this setting the Sally module
has either dimension $d$ or $0$.

\begin{Proposition}\label{dimension}
Let $(R, {\mathfrak m})$ be a local Noetherian ring of dimension
$d>0$ with infinite residue field and let $J$ be a minimal
reduction of ${\mathfrak m}$. Then the Sally module
$S_J({\mathfrak m})$ of ${\mathfrak m}$ with respect to $J$ has
dimension $d$ if and only if $e_1({\mathfrak m})-e_0({\mathfrak
m})-e_1(J)+1$ is strictly positive. In this event, the
multiplicity of $S_J({\mathfrak m})$ is exactly $e_1({\mathfrak
m})-e_0({\mathfrak m})-e_1(J)+1$.
\end{Proposition}
\demo We compute the Hilbert function/polynomial of the Sally
module $S_J({\mathfrak m})$. Chasing lengths in the short exact
sequences
\[
0 \rar {\mathfrak m}J^{n-1}/J^n \lar {\mathfrak m}^n/J^n \lar
\fbox{${\mathfrak m}^n/{\mathfrak m}J^{n-1}$} \rar 0,
\]
\[
0 \rar {\mathfrak m}J^{n-1}/J^n \lar J^{n-1}/J^n \lar
J^{n-1}/{\mathfrak m}J^{n-1} \rar 0
\]
leads to the following equality that provides a formula for the
length of the component of degree $n-1$ of the Sally module
$S_J({\mathfrak m})$
\[
\lambda({\mathfrak m}^n/{\mathfrak m}J^{n-1}) =
\lambda(R/J^{n-1})-\lambda(R/{\mathfrak
m}^n)+\lambda(J^{n-1}/{\mathfrak m}J^{n-1}).
\]
Now observe that, for $n \gg 0$, both $\lambda(R/{\mathfrak m}^n)$
and $\lambda(R/J^{n-1})$ can be replaced with their respective
Hilbert-Samuel polynomials
\begin{eqnarray*}
\lambda(R/{\mathfrak m}^n) \!\!\!\! & = & \!\!\!\! e_0({\mathfrak
m}) {{n+d-1}\choose d} - e_1({\mathfrak m}){{n+d-2}\choose d-1} +
\cdots + (-1)^d e_d({\mathfrak m}), \\
\lambda(R/J^{n-1}) \!\!\!\! & = & \!\!\!\! e_0(J) {{n+d-2}\choose
d} - e_1(J) {{n+d-3}\choose d-1} + \cdots + (-1)^d e_d(J).
\end{eqnarray*}
On the other hand $J^{n-1}/{\mathfrak m}J^{n-1}$ is the component
of degree $n-1$ of the special fiber ring ${\mathcal F}(J)$ of
$J$, which is a polynomial ring in $d$ variables with coefficients
over the residue field. In particular we have that
\[
\lambda(J^{n-1}/{\mathfrak m}J^{n-1}) = {n+d-2 \choose d-1}.
\]
Finally, using the fact that $e_0({\mathfrak m})=e_0(J)$ as shown
in \cite[Section 1, Theorem 1]{NR}, since $J$ is a reduction of
${\mathfrak m}$, and the combinatorial identity ${p \choose q} +
{p \choose q+1} = {p+1 \choose q+1}$ for non-negative integers $p$
and $q$, we have
\[
\lambda({\mathfrak m}^n/{\mathfrak m}J^{n-1}) = s_0 \, {n+d-2
\choose d-1} -s_1 {n+d-3 \choose d-2} +\cdots +(-1)^{d-1}s_{d-1},
\]
where \ $s_0=e_1({\mathfrak m})-e_0({\mathfrak m})-e_1(J)+1$ \ and
\ $s_i= e_{i+1}({\mathfrak m})-e_i(J)-e_{i+1}(J)$ \ for $i=1,
\ldots d-1$. This proves that the dimension of $S_J({\mathfrak
m})$ is $d$ if and only if $e_1({\mathfrak m})-e_0({\mathfrak
m})-e_1(J)+1$ is strictly positive. \QED

\bigskip

In a local Cohen-Macaulay setting the Sally module has the same
dimension as the ambient ring, unless it is the trivial module.
This is no longer the case in an arbitrary local Noetherian ring
of positive dimension, as the next example taken from
\cite[4.2]{GN} shows.

\begin{Example}\label{ex-1}
{\rm Let $k$ be a field and let $S=k[X,Y,Z,W]$ be the polynomial
ring in $4$ variables over $k$. Define $T=S/(X^2,Y)\cap(Z,W)$,
$M=T_{+}$ and $R=T_M$ and ${\mathfrak m}=MR$. Let $x$, $y$, $z$
and $w$ denote the images in $R$ of $X$, $Y$, $Z$ and $W$,
respectively. The ring $R$ is a two-dimensional local ring such
that the Sally module $S_J({\mathfrak m})$ of ${\mathfrak m}$ with
respect to $J=(x-z, y-w)$ has dimension one. In fact
$e_1({\mathfrak m})-e_0({\mathfrak m})-e_1(J)+1=0$ whereas
$-s_1=-e_2({\mathfrak m})+e_1(J)+e_2(J)=1$. }
\end{Example}

\smallskip

In Theorem~\ref{E-V} below we use the techniques of \cite{CPV2} to
generalize to a local Noetherian setting an estimate of J. Elias
and G. Valla \cite[Section 2]{EV}, which involves the Hilbert
coefficients of ${\mathfrak m}$ and the embedding codimension of
$R$. Also, in Example~\ref{Macaulay} we use a well-known example,
even studied by F.S. Macaulay as early as 1916, to provide an
instance that illustrates when equality in Theorem~\ref{E-V} is
attained.

\begin{Theorem}\label{E-V}
Let $(R, {\mathfrak m})$ be a local Noetherian ring of dimension
$d>0$ with infinite residue field. Then \ $2e_0({\mathfrak
m})-e_1({\mathfrak m})+e_1(J) \leq \mu({\mathfrak m})-d+2$, \
where $J$ is any minimal reduction of ${\mathfrak m}$.
\end{Theorem}
\demo Let $J$ be a minimal reduction of ${\mathfrak m}$ and write
${\mathfrak m}=(J, a_1, \ldots, a_{m-d})$, where $m$ denotes the
minimal number of generators of ${\mathfrak m}$. We now consider
the Sally module $S_J({\mathfrak m})$ of $I$ with respect to $J$
defined by means of the exact sequence introduced in \cite[proof
of 2.1]{CPV2}
\[
{\mathcal R}(J) \oplus  {\mathcal R}(J)^{m-d}[-1]
\stackrel{\varphi}{\longrightarrow} {\mathcal R}({\mathfrak m})
\longrightarrow S_J({\mathfrak m})[-1] \rightarrow 0,
\]
where $\varphi$ is the map defined by $\varphi(r_0, r_1, \ldots,
r_{m-d}) = r_0+r_1a_1t+\cdots+r_{m-d}a_{m-d}t$, for any element
$(r_0, r_1, \ldots, r_{m-d})\in {\mathcal R}(J) \oplus {\mathcal
R}(J)^{m-d}[-1]$. Tensoring the above exact sequence with
$R/{\mathfrak m}$ yields the bottom row in the diagram
\[
\begin{array}{lcr} \vspace{0.1cm}
& S_J({\mathfrak m})[-1] & \\ \vspace{0.1cm}
& \downarrow & \\
\vspace{0.1cm} {\mathcal F}(J) \oplus {\mathcal F}(J)^{m-d}[-1]
\longrightarrow {\mathcal G}({\mathfrak m}) \longrightarrow &
\!\!\!\! S_J({\mathfrak
m})[-1] \otimes R/{\mathfrak m} & \!\!\!\! \rightarrow 0. \\
& \downarrow & \\
& 0 &
\end{array}
\]
Hence we obtain the following multiplicity $($degree$)$ estimate
\[
e_0({\mathfrak m}) \leq e_1({\mathfrak m})-e_0({\mathfrak
m})-e_1(J)+1 + \deg ({\mathcal F}(J) \oplus {\mathcal
F}(J)^{m-d}[-1]).
\]
Notice that the contribution involving $S_J({\mathfrak m})$ only
occurs if its dimension is $d$, in which case we use the result of
Proposition~\ref{dimension}. On the other hand, ${\mathcal F}(J)
\oplus {\mathcal F}(J)^{m-d}[-1]$ is a free ${\mathcal
F}(J)$-module of rank $m-d+1$. Thus, its multiplicity is $m-d+1$,
since ${\mathcal F}(J)$ is isomorphic to a polynomial ring. The
result now easily follows. \QED

\bigskip

In the rest of the section we restrict our attention to the case
of a local Buchsbaum ring. We briefly review some basic notions,
whereas we refer the reader to the monograph of J. St\"uckrad and
W. Vogel for a comprehensive treatment of the subject \cite{SV}.
In short, the theory of Buchsbaum rings is a natural
generalization of the concept of a Cohen-Macaulay ring and started
in a remarkable series of papers by J. Stuckrad and W. Vogel to
answer negatively a problem of D.A. Buchsbaum. A local Noetherian
ring $(R, {\mathfrak m})$ of positive dimension $d$ is said to be
a {\it Buchsbaum ring} if and only if there exists a non-negative
integer $I(R)$ such that $\lambda(R/J) - e_0(J) = I(R)$ for every
system of parameters $J=(x_1,\dots, x_d)$ of $R$. In this setting
any such $J$ is no longer generated by a regular sequence $($as in
the Cohen-Macaulay case, where $I(R)=0)$, but by a {\it
$d$-sequence}.
The number $I(R)$ is the so-called {\it Buchsbaum invariant} of
$R$ and has an explicit description either in terms of the lengths
of the local cohomology modules $H_{\mathfrak m}^i(R)$, for $i=0,
\ldots, d-1$ $($as they are annihilated by the maximal ideal
${\mathfrak m})$, or in terms of the higher Hilbert coefficients
of any system of parameters $J$ of $R$. Namely, we have that
\[
I(R) = \sum_{i=0}^{d-1} {d-1\choose i} \lambda(H_{\mathfrak
m}^i(R)) \qquad {\rm or} \qquad I(R) = \sum_{i=1}^d (-1)^i e_i(J),
\]
for every parameter ideal $J$ of $R$. Even more surprisingly, in a
local Buchsbaum ring $R$ one has that the Hilbert coefficients
$e_i(J)$, for $i=1, \ldots, d$, do not depend on the system of
parameters $J$ but only on the ring $R$. $($Note that, in a local
Cohen-Maculay ring, $e_i(J)=0$ for $i>0.)$ In particular, we
observe that if $R$ is a Buchsbaum ring then the formula in
Theorem~\ref{E-V} does not depend on the reduction $J$ of
${\mathfrak m}$, but solely on the ring $R$. Indeed, for any
parameter ideal $J$ of a local Buchsbaum ring $R$ one has that
\[
-e_1(J) = \sum_{i=0}^{d-1} {d-2 \choose i-1} \lambda(H_{\mathfrak
m}^i(R)),
\]
where ${p \choose -1}$ is either 0 when $p \not= -1$, or 1 when
$p=-1$ $($see \cite[2.7(ii)]{SV}$)$.

\begin{Example}\label{Macaulay} {\rm
The curve $X$ in ${\mathbb P}^3$ given parametrically by
$\{s^4,s^3u,su^3,u^4\}$ is such that the local ring of the affine
cone over $X$ at the vertex is a Buchsbaum ring of dimension $2$
with invariant $I(R)=1=-e_1(J)$ for every system of parameters $J$
of $R$. One can verify that equality holds in the bound
established in Corollary~\ref{E-V}, as $e_0({\mathfrak m})=4$ and
$e_1({\mathfrak m})=3$. Moreover, the maximal ideal ${\mathfrak
m}$ has reduction number $2$, the associated graded ring
${\mathcal G}({\mathfrak m}) \cong R$ is Buchsbaum, the Sally
module $S_J({\mathfrak m})$ of ${\mathfrak m}$ with respect to a
minimal reduction $J$ is two-dimensional with multiplicity $1$ and
${\mathfrak m}\,S_J({\mathfrak m})=0$. }
\end{Example}

The above example and other similar ones led us to ask the
following question, which is motivated by an analogous result due
to J. Elias and G. Valla in a local Cohen-Macaulay setting
\cite[2.1]{EV}.

\begin{Question}
{\rm Let $(R, {\mathfrak m})$ be a local Buchsbaum ring and
suppose that \ $2e_0({\mathfrak m})-e_1({\mathfrak m})+e_1(J) =
\mu({\mathfrak m})-d+2$, \ where $J$ is any minimal reduction of
${\mathfrak m}$. Is the associated graded ring ${\mathcal
G}({\mathfrak m})$ always Buchsbaum? }
\end{Question}

We end this section by showing that the dimension of the Sally
module is rather extremal $($either $d$ or $0)$ in the case in
which the ambient ring is Buchsbaum. Our proof uses a remarkable
generalization due to S. Goto and K. Nishida \cite[1.1]{GN} of the
result by C. Huneke \cite[2.1]{H2} and A. Ooishi \cite[3.2 and
3.3]{O} quoted in the introduction.

\begin{Proposition}\label{dimbuch}
Let $(R, {\mathfrak m})$ be a local Buchsbaum ring of dimension
$d>0$ with infinite residue field and let $J$ be a minimal
reduction of ${\mathfrak m}$. Then the Sally module
$S_J({\mathfrak m})$ of ${\mathfrak m}$ with respect to $J$ has
either dimension $d$ or $0$.
\end{Proposition}
\demo As shown in the proof of Proposition~\ref{dimension}, if
$e_1({\mathfrak m})-e_0({\mathfrak m})-e_1(J)+1>0$ then the Sally
module $S_J({\mathfrak m})$ has dimension $d$. Suppose now that
$e_1({\mathfrak m})-e_0({\mathfrak m})-e_1(J)+1=0$. Let $H$ denote
the zero-th local cohomology module $H^0_{\mathfrak m}(R)$ of $R$
with support in ${\mathfrak m}$ and let $\,{}^{{}^{\tratto}}$
denote images in the ring $\overline{R}=R/H$. The two Sally
modules $S_J({\mathfrak m})$ and
$S_{\overline{J}}(\overline{\mathfrak m})$ are related by the
short exact sequence
\[
0 \rar K=\bigoplus_{n\geq 2}^{\infty} \frac{H \cap {\mathfrak m}^n
+{\mathfrak m}J^{n-1}}{{\mathfrak m}J^{n-1}} \lar S_J({\mathfrak
m}) \lar S_{\overline{J}}(\overline{\mathfrak m}) \rar 0.
\]
Since $H$ is Artinian and therefore $H \cap {\mathfrak m}^n=0$ for
$n$ sufficiently large, it follows that the kernel $K$ has only
finitely many components, hence it is Artinian. On the other hand,
by \cite[1.1]{GN} we have that ${\mathfrak m}^2 \subset
J{\mathfrak m}+H$. Hence $S_{\overline{J}}(\overline{\mathfrak
m})=0$ so that $S_J({\mathfrak m}) = K$ is Artinian. \QED

\bigskip

We also observe that in the case of a zero-dimensional Sally
module all the Hilbert coefficients of ${\mathfrak m}$ only depend
on the multiplicity of ${\mathfrak m}$ and the local Buchsbaum
ring $R$.

\begin{Corollary}
Let $(R, {\mathfrak m})$ be a local Buchsbaum ring of dimension
$d>0$ with infinite residue field and let $J$ be a minimal
reduction of ${\mathfrak m}$. If $e_1({\mathfrak
m})=e_0({\mathfrak m})+e_1(J)-1$ then
\[
e_i({\mathfrak m}) = e_{i-1}(J)+e_i(J),
\]
for $i=2, \ldots, d$.
\end{Corollary}

\smallskip

\section{\bf Sally modules of ${\mathfrak m}$-primary ideals and
special fiber rings \hfill\break}

It is natural to ask what happens in the case of an arbitrary
${\mathfrak m}$-primary ideal $I$ of a local Noetherian ring $R$.
We address this issue next. In this case we can only give an upper
bound for the multiplicity of the Sally module $S_J(I)$ of $I$
with respect to a minimal reduction $J$, unless the ideal $I$
contains the zero-th local cohomology module $H_{\mathfrak
m}^0(R)$ of $R$ with support in ${\mathfrak m}$. Nevertheless,
this enables us to obtain general multiplicity estimates for the
special fiber ring ${\mathcal F}(I)$: see Theorem~\ref{fiber} and
Theorem~\ref{fiber2}. In a Cohen-Macaulay setting, the estimates
had been previously obtained in a paper by W.V. Vasconcelos
\cite[2.4]{mrn} and in a joint work with C. Polini, W.V.
Vasconcelos \cite[2.1, 2.2]{CPV2}.

\begin{Proposition}\label{sally}
Let $(R, {\mathfrak m})$ be a local Noetherian ring of dimension
$d>0$ with infinite residue field and let $I$ be an ${\mathfrak
m}$-primary ideal. Then a $d$-dimensional Sally module $S_J(I)$ of
$I$ with respect to a minimal reduction $J$ has multiplicity at
most \ $e_1(I)-e_0(I)-e_1(J)+\lambda(R/I)$, \ with equality if and
only if $I$ contains $H_{\mathfrak m}^0(R)$.
\end{Proposition}
\demo In a similar fashion as in the proof of
Proposition~\ref{dimension}, we obtain a formula for the length of
the component of degree $n-1$ of the Sally module $S_J(I)$
\[
\lambda(I^n/IJ^{n-1}) =
\underbrace{(e_1(I)-e_0(I)-e_1(J)+\hat{e}_0(J,I))}_{\text{leading
coefficient}} \, {n+d-2 \choose d-1} + \text{ lower terms},
\]
where $\hat{e}_0(J,I)$ is the multiplicity of graded module
${\mathcal G}(J) \otimes R/I$. Observe that the statement follows
once we show that $\hat{e}_0(J,I) \leq \lambda(R/I)$, with
equality if and only if $I$ contains $H_{\mathfrak m}^0(R)$. We
prove this claim by induction on the dimension $d$ of $R$. If
$d=1$, let $H$ denote the zero-th local cohomology module
$H^0_{\mathfrak m}(R)$ of $R$ with support in ${\mathfrak m}$ and
let $\,{}^{{}^{\tratto}}$ denote images in the one-dimensional
local Cohen-Macaulay ring $\overline{R}=R/H$. Since $H \cap J^n$
is eventually the zero ideal as $H$ is Artinian, it is not
difficult to verify that $J^n/IJ^n \cong
\overline{J}^n/\overline{I}\overline{J}^{n}$ for $n \gg 0$. Thus
$\hat{e}_0(J,I) =
\hat{e}_0(\overline{J},\overline{I})=\lambda(\overline{R}/\overline{I})
\leq \lambda(R/I)$, with equality if and only if $I$ contains $H$.

Suppose now $d>1$. We can find an element $x$ which belongs to $J
\setminus {\mathfrak m}J$ and whose image is superficial in
${\mathcal G}(J)\otimes R/I$. Let $'$ denote images in the ring
$R'=R/(x)$, which has dimension $d-1$. Observe that
$\hat{e}_0(J,I) = \hat{e}_0(J',I')$, as the image of $x$ is
superficial in ${\mathcal G}(J)\otimes R/I$ $($see
\cite[22.6]{Na}$)$. Thus, by inductive hypothesis we conclude that
$\hat{e}_0(J,I) = \hat{e}_0(J',I') \leq
\lambda(R'/I')=\lambda(R/I)$. On the other hand, by induction we
have that $\hat{e}_0(J',I')=\lambda(R'/I')$ if and only if
$H^0_{\mathfrak m}(R') \subset I'$. Our assertion now follows as
$H^0_{\mathfrak m}(R)R' \subset H^0_{\mathfrak m}(R')$.  \QED

\medskip

\begin{Remark}{\rm
If in addition to the assumptions in Proposition~\ref{sally} all
the local cohomology modules $H_{\mathfrak m}^i(R)$ have finite
length for $i < d$, S. Goto and K. Nishida have shown that
\[
-e_1(J) \leq \sum_{i=0}^{d-1} {d-2 \choose i-1}
\lambda(H_{\mathfrak m}^i(R)).
\]
$($see \cite[2.4]{GN}$)$. This yields another bound on the
multiplicity of a $d$-dimensional Sally module which depends
solely on the ideal and still reduces to the classical bound of
W.V. Vasconcelos in the Cohen-Macaulay case. }
\end{Remark}

\medskip

The next result is similar to the one in
Proposition~\ref{dimbuch}.

\begin{Proposition}\label{dimbuch-primary}
Let $(R, {\mathfrak m})$ be a local Buchsbaum ring of dimension
$d>0$ with infinite residue field and let $I$ be an ${\mathfrak
m}$-primary ideal containing $H_{\mathfrak m}^0(R)$. Then the
Sally module $S_J(I)$ of $I$ with respect to a minimal reduction
$J$ has either dimension $d$ or $0$.
\end{Proposition}
\demo The proof is similar to the one of
Proposition~\ref{dimbuch}. In fact, if $S_J(I)$ has not dimension
$d$ then then $e_1(I)-e_0(I)-e_1(J)+\lambda(R/I)=0$ by
Proposition~\ref{sally} and hence $I^2 \subset JI + H_{\mathfrak
m}^0(R)$ by \cite[1.1]{GN}. \QED

\medskip

We now use the previous results to obtain estimates on the
multiplicity of the special fiber ring of any ${\mathfrak
m}$-primary ideal.

\begin{Theorem}\label{fiber} Let $(R, {\mathfrak m})$ be a local
Noetherian ring of dimension $d>0$ with infinite residue field and
let $I$ be an ${\mathfrak m}$-primary ideal. Then the multiplicity
$f_0(I)$ of the special fiber ring ${\mathcal F}(I)$ of $I$ is at
most \ $e_1(I)-e_0(I)-e_1(J)+\lambda(R/I)+\mu(I)-d+1$, \ where $J$
is any minimal reduction of $I$.
\end{Theorem}
\demo Let $J$ be any minimal reduction of $I$ and let $m$ denote
the minimal number of generators of $I$. We proceed in a similar
fashion as in the proof of Theorem~\ref{E-V}. Namely, tensoring
with $R/{\mathfrak m}$ the defining sequence of the Sally module
$S_J(I)$ of $I$ with respect to $J$
\[
{\mathcal R}(J) \oplus  {\mathcal R}(J)^{m-d}[-1]
\stackrel{\varphi}{\longrightarrow} {\mathcal R}(I)
\longrightarrow S_J(I)[-1] \rightarrow 0
\]
yields the following exact sequence
\[
{\mathcal F}(J) \oplus {\mathcal F}(J)^{m-d}[-1] \longrightarrow
{\mathcal F}(I) \longrightarrow S_J(I)[-1] \otimes R/{\mathfrak m}
\rightarrow 0.
\]
Our assertion now follows after taking into account the estimate
of Proposition~\ref{sally}. \QED

\medskip

\begin{Theorem}\label{fiber2}
Let $(R, {\mathfrak m})$ be a local Buchsbaum ring of dimension
$d>0$ with infinite residue field and let $I$ be an ${\mathfrak
m}$-primary ideal. Then the multiplicity $f_0(I)$ of the special
fiber ring ${\mathcal F}(I)$ of $I$ is at most \ $e_1(I) + I(R)
-e_1(J)+ 1$, \ where $J$ is any minimal reduction of $I$. \
Furthermore, if the bound is attained then the ideal $I$ has
minimal multiplicity in the sense of\/ {\rm S. Goto}, that is
${\mathfrak m}I={\mathfrak m}J$ for any minimal reduction $J$ of
$I$.
\end{Theorem}
\demo Our assertions follow from Theorem~\ref{fiber} and the
equality due to K. Yamagishi $e_0(I) = \mu(I)-d +\lambda(R/I) -
I(R)+\lambda({\mathfrak m}I/{\mathfrak m}J)$ \cite[2.5, 2.6]{Y2},
which generalizes the classical result of S. Abhyankar. \QED

\medskip

The next example shows that the ideal $I$ may have minimal
multiplicity even if the inequality in Theorem~\ref{fiber2} is
strict.

\begin{Example}{\rm
Let $k$ be a field and let $S=k[X_1, X_2, X_3, X_4, V, A_1, A_2,
A_3]$ be the polynomial ring in $8$ variables over $k$. and put
\[
{\mathfrak a} = (X_1, X_2, X_3)^2+(X_4^2)+(X_1V, X_2V, X_3V,
X_4V)+(V^2-A_1X_1-A_2X_2-A_3X_3).
\]
Define $T=S/{\mathfrak a}$, $M=T_{+}$ and $R=T_M$ and ${\mathfrak
m}=MR$. Let $x_i$, $v$ and $a_j$ denote the images of $X_i$, $V$
and $A_j$ modulo ${\mathfrak a}$, respectively. It follows from
\cite[Section 4]{GS2} that $R$ is a local Buchsbaum ring of
dimension $3$ with $I(R)=1$. Moreover, the ideal $I=J \colon
{\mathfrak m}$, where $J=(a_1, a_2, a_3)R$, is such that
$f_0(I)=5$, $e_1(I)=3$ and $e_1(J)=-1$. Thus we have a strict
inequality in the bound of Theorem~\ref{fiber2}. However, the
ideal $I$ satisfies $I^3=JI^2$ and has minimal multiplicity, that
is ${\mathfrak m}I={\mathfrak m}J$. }
\end{Example}

We conclude with a generalization of \cite[2.9]{CPV2} to a non
Cohen-Macaulay setting.

\begin{Proposition}\label{E-V2}
Let $(R, {\mathfrak m})$ be a local Noetherian ring of dimension
$d>0$ with infinite residue field and let $I$ be an ${\mathfrak
m}$-primary ideal. Then
\[
2e_0(I)-e_1(I)+e_1(J) \leq \lambda(R/I)\left(\mu(I)-d+2 \right),
\]
where $J$ is any minimal reduction of $I$.
\end{Proposition}
\demo Repeat the proof of Theorem~\ref{E-V}, using the estimate
given in Proposition~\ref{sally}. \QED

\end{document}